\documentclass[]{interact}
\usepackage{epstopdf}
\usepackage[caption=false]{subfig}

\usepackage[natbibapa,nodoi]{apacite}
\theoremstyle{plain}

\theoremstyle{definition}

\theoremstyle{remark}

\usepackage{amsmath}

\usepackage{tikz}
\usepackage{color}
\usepackage{array}
\usepackage{multirow}
\usepackage{float}
\usepackage{makecell}
\usepackage{fancyhdr}
\usepackage{comment}
\usepackage{calrsfs}
\DeclareMathAlphabet{\pazocal}{OMS}{zplm}{m}{n}

\usepackage{listings}
\usepackage{calrsfs} 
\usepackage{natbib}
\usepackage{caption}
\usepackage{subfig}
\usepackage{algorithm}
\usepackage{algpseudocode}
\usepackage{etoolbox}

\makeatletter
\def\BState{\State\hskip-\ALG@thistlm}

\newcommand{\bzero}{ \mathbf{0} }

\newcommand{\blambda}{ \boldsymbol{\lambda} }

\usepackage{calrsfs}

\newcommand{\bx}{ \mathbf{x} }

\newcommand{\bc}{ \mathbf{c} }
\newcommand{\bu}{ \mathbf{u} }

\newcommand{\bA}{ \mathbf{A} }
\newcommand{\bb}{ \mathbf{b} }
\usepackage[english]{babel} 
\usepackage{amsmath,amsfonts,amsthm,bm} 
\DeclareMathOperator*{\argmin}{arg\,min}
\usepackage{graphicx}
\usepackage{amssymb}

\usepackage{lineno}
\usepackage{graphicx}
\usepackage{natbib}
\usepackage{url}
\usepackage{amstext}
\usepackage{amssymb}
\usepackage{tikz}
\allowdisplaybreaks
\usepackage{mathtools}

\newcommand*{\algrule}[1][\algorithmicindent]{%
  \makebox[#1][l]{%
    \hspace*{.2em}
    \vrule height .75\baselineskip depth .25\baselineskip
  }
}

\newcount\ALG@printindent@tempcnta
\def\ALG@printindent{%
    \ifnum \theALG@nested>0
    \ifx\ALG@text\ALG@x@notext
    \else
    \unskip
    \ALG@printindent@tempcnta=1
    \loop
    \algrule[\csname ALG@ind@\the\ALG@printindent@tempcnta\endcsname]%
    \advance \ALG@printindent@tempcnta 1
    \ifnum \ALG@printindent@tempcnta<\numexpr\theALG@nested+1\relax
    \repeat
    \fi
    \fi
}
\patchcmd{\ALG@doentity}{\noindent\hskip\ALG@tlm}{\ALG@printindent}{}{\errmessage{failed to patch}}
\patchcmd{\ALG@doentity}{\item[]\nointerlineskip}{}{}{} 

\algnewcommand{\algorithmicforeach}{\textbf{for each}}
\algdef{SE}[FOR]{ForEach}{EndForEach}[1]
  {\algorithmicforeach\ #1\ \algorithmicdo}
  {\algorithmicend\ \algorithmicforeach}

\begin{document}


\title{Comparing Inverse Optimization and Machine Learning Methods for Imputing a Convex Objective Function}

\author{
\name{Elaheh H.Iraj\thanks{Email address: elaheh.hosseiniiraj@concordia.ca, daria.terekhov@concordia.ca (corresponding author)} and Daria Terekhov}
\affil{Department of Mechanical, Industrial and Aerospace Engineering, Gina Cody School of Engineering and Computer Science, Concordia University, Montr\'eal, QC,
Canada}
}

\maketitle

\begin{abstract}
Inverse optimization (IO) aims to determine optimization model parameters from observed decisions. However, IO is not part of a data scientist's toolkit in practice, especially as many general-purpose machine learning packages are widely available as an alternative. When encountering IO, practitioners face the question of when, or even whether, investing in developing IO methods is worthwhile. Our paper provides a starting point toward answering these questions, focusing on the problem of imputing the objective function of a parametric convex optimization problem. We compare the predictive performance of three standard supervised machine learning (ML) algorithms (random forest, support vector regression and Gaussian process regression) to the performance of the IO model of \citet{Keshavarz11}. While the IO literature focuses on the development of methods tailored to particular problem classes, our goal is to evaluate general `out-of-the-box' approaches. Our experiments demonstrate that determining whether to use an ML or IO approach requires considering 
(i) the training set size, (ii) the dependence of the optimization problem on external parameters, (iii) the level of confidence with regards to the correctness of the  optimization prior, and (iv) the number of critical regions in the solution space.
\end{abstract}

\begin{keywords}
inverse optimization; machine learning; data efficiency; parametric optimization; multi-parametric programming
\end{keywords}

\section{Introduction}\label{sec:intro}
Inverse optimization (IO) imputes missing optimization model parameters from data that represents minimally sub-optimal solutions of that unknown optimization problem (e.g., past decisions of an optimizing agent). In the academic literature, IO has been shown successful in a variety of problems, including medical decision making \citep{chan2014generalized}, electricity demand forecasting \citep{Gallego}, the household activity pattern problem \citep{chow2012inverse} and economic lot-sizing \citep{egri2014inverse}. However, IO is rarely used by practitioners. 

Previous work has established an analogy between IO and regression \citep{GIO}, and provided a statistical inference perspective on IO \citep{aswani2018inverse}. More generally, one can notice that \emph{imputing parameter values from data} can be viewed as a machine learning problem \citep{deepIO}, a perspective that we adopt and explore in this paper. Viewing an inverse (parametric) optimization problem as a learning problem raises questions that have previously not been explored in the literature and that are of practical interest. First, how well would classic machine learning methods perform on such problems? Second, what characteristics would make a problem challenging for classic machine learning methods and would require the more specialized methods of IO? Our interactions with data science teams of two large companies have led us to believe that answering these questions is essential if a data scientist wants to consider adding IO to their toolkit.

While most of the IO literature develops specialized approaches, our goal is to compare fairly general methods, as we believe these are more likely to be used by practitioners. After establishing our perspective on inverse parametric optimization as a learning problem, we experimentally compare the performance of a classic IO method with three out-of-the-box machine learning (ML) methods, namely Gaussian process regression, random forest and support vector regression. We use the IO method of \citet{Keshavarz11} as it is well suited and easy to apply to inverse parametric convex optimization problems.  We show that the chosen ML methods can indeed perform well on a problem that has been well-studied in the IO literature. To identify the characteristics that would make an IO problem challenging for ML  (highlighting the need for sophisticated IO models), we conduct experiments on random parametric optimization problems (POPs) generated using the parametric optimization toolbox of~\citet{POP}. To the best of our knowledge, no previous papers have compared IO and ML methods on this class of problems. Our experiments with randomly-generated POPs demonstrate that the choice of an ML or IO approach should depend on (i) the size of the training set, (ii) the nature of the dependence of the optimization problem on external parameters, (iii) the level of confidence with regards to the correctness optimization prior, (iv) the number of critical regions in the solution space of the POP. 


\section{Background and Related Work}\label{sec:lit_review}

To define an IO problem, we first need to postulate a \emph{forward} optimization problem, i.e., the problem whose solutions, perhaps corrupted by noise, we observe. The initial focus of the IO literature was to impute coefficients of forward optimization problems that do not have any dependence on an external parameter (non-parametric). \citet{burton1992instance, burton_correlated94} study the inverse shortest path problem in which the goal is to find a minimal perturbation of the arc costs to make an observed set of paths optimal. Given an observed value for decision variables $\bx$, \citet{Ahuja01} find a $\bc$ vector of minimal distance to a known vector $\hat{\bc}$ that would make $\bx$ optimal. \citet{chan2014generalized} and \citet{GIO} study the inverse linear programming problem when the observed data is noisy and no prior $\hat{\bc}$ is given. \citet{tavasliouglu2018}  characterize the inverse feasible region: the set of objectives that would make a given set of feasible solutions to a linear program optimal. \citet{shahmoradi2020} introduce a way to measure how sensitive the set of optimal cost vectors is to changes in a given data set.

A parametric optimization problem (also known as a multi-parametric programming problem) is a family of optimization problems parametrized by an independent parameter; the forward parametric problem involves finding a function that would map values of the independent parameter to optimal solutions \citep{pistikopoulos2011multi}. As an example, consider the following parameteric linear program with independent parameters $(u_1, u_2)$: $\min (c_{1}+c_{2}u_{1})x_{1}+(c_{3}+c_{4}u_{2})x_{2} \textrm{ s.t. } a_1x_1 + a_2x_2 \ge b_1$, 
where instantiating $u_1$ and $u_2$ leads to a standard (non-parametric) optimization problem. There is substantial interest in the inverse parametric optimization problem  \citep{Keshavarz11,chow2012inverse,aswani2018inverse,esfahani2018data,Gallego,kovacs2019parameter,deepIO,tan2020}. In such a problem, pairs $(\bu_k, \bx_k)$ where $\bu_k$ is an independent parameter for observation $k$ and $\bx_k$ is an observed value for the decision variables $\bx$ are given; the goal is to impute objective function and/or constraint coefficients that would make the observed points optimal or near-optimal solutions of the given parametric forward optimization problem. 

Classic IO methods \citep{burton1992instance,Ahuja01,Keshavarz11} rely on some form of optimality conditions, such as the Karush-Kuhn-Tucker (KKT) conditions. More recently, methods based on a combination of ML and IO have emerged \citep{deepIO,tan2020,babier2019}. \citet{aswani2019behavioral} and \citet{fernandez2019ev} include comparison of an IO approach to several ML algorithms in the context of a particular application. Related work in the ML literature includes, most notably, structured prediction \citep{taskar2005learning}. Despite this related ML-based work, a comparison of classic methods on general parametric IO problems has not previously been done and the questions posed in the introduction have not been answered.

\section{Problem Definition}\label{sec:problem_def}
We study IO for \textit{parametric optimization problems} (POPs), which are `compatible' with an ML perspective unlike the non-parametric optimization ones. POPs are mathematical programs where the objective function and constraints are functions of one or multiple external parameters \citep{pistikopoulos2011multi,POP}.  
The \emph{forward} parametric optimization problem (FOP) of interest to us is: 
\begin{equation}\label{eq:FOP1}
\begin{alignedat}{1}
\mathbf{FOP} (\textbf{u}, \textbf{c}) : \underset{\textbf{x}}{\text{minimize}} & \quad  f(\textbf{u}, \textbf{x}, \textbf{c})  \\
\text{s.t.}
 & \quad  g(\textbf{u},\textbf{x})\leq \textbf{0}, \\
\end{alignedat}
\end{equation}
\noindent where $\bx \in \mathbb{R}^n$ are the decision variables, $\bu \in \mathbb{R}^v$ are the independent parameters (`features' in ML parlance), and $\bc \in \mathbb{R}^n$ are the coefficients of the objective function $f$. We assume $f$ and $g$ are differentiable and convex in $\textbf{x}$ for each value of $\textbf{u}$. We let $\pazocal{J}$ be the index set of variables with $|\pazocal{J}|=n$, and $\pazocal{I}$ be the index set of constraints with $|\pazocal{I}|=m$; we let $\bzero$ denote the column vector of zeros of the required dimension.   

In an \emph{inverse} parametric optimization problem, observations of pairs $\pazocal{D} = \{ (\hat{\textbf{x}}_{1}, \hat{\bu}_1), \dots, (\hat{\textbf{x}}_{K}, \hat{\bu}_K)\}$ are given; the goal is to find a $\bc$ that would make $\hat{\textbf{x}}_{k}$ minimally sub-optimal for the corresponding $\mathbf{FOP} (\hat{\bu}_k, \textbf{c})$, for each
$k \in \pazocal{K} = \{1, \dots, K\}$:

\begin{equation} \label{eq:IO2}
\underset{\textbf{c}}{\text{{minimize}}} \quad \Big\{ \sum_{k \in \pazocal{K}} \pazocal{L}(\hat{\mathbf{x}}_k, \mathbf{x}_k^{pred}) \enskip | \enskip \textbf{x}_{k}^{pred} \in \underset{\textbf{x} \in \pazocal{X}_{\hat{\textbf{u}}_{k}}}{\argmin} \enskip f(\hat{\textbf{u}}_{k}, \textbf{x}, \textbf{c}) \Big\} 
\end{equation}
\noindent where $\pazocal{L}$ is any given loss function, $\mathbf{x}^{pred}_k$ is an optimal solution to the fitted forward optimization model instantiated at $\mathbf{u} = \hat{\mathbf{u}}_k$, and $\pazocal{X}_{\hat{{\bu}}_k} = \{\mathbf{x} \; | \; g(\hat{\mathbf{u}}_k,\mathbf{x}) \le \mathbf{0}\}$ is the feasible region corresponding to $\hat{\bu}_k$. There are various ways to reformulate \eqref{eq:IO2}; the particular model we use here is the model proposed by \citet{Keshavarz11}, which formulates the optimality criterion via KKT optimality conditions (see appendix). 

The goal of IO could be explanatory, i.e., we aim to characterize the process that generated the data, or predictive, i.e., we want to fit a model that will predict the actions of an optimizing agent/process. That is, having imputed $\bc$ and given $\textbf{u}^{new}$ we could use $\mathbf{FOP}(\textbf{u}^{new}, \textbf{c})$ to predict $\textbf{x}^{new}$. Fitting a model to data can also result in better prescriptive performance, i.e., if instead of predicting what an optimizing agent/process would do, we want to prescribe an optimal decision under new circumstances $\textbf{u}^{new}$ (this latter perspective being more consistent with classical optimization). 

\section{Inverse Optimization Problem as a Learning Problem}\label{sec:IOasML}
\begin{figure}[H]
\centering
\includegraphics[scale=0.14]{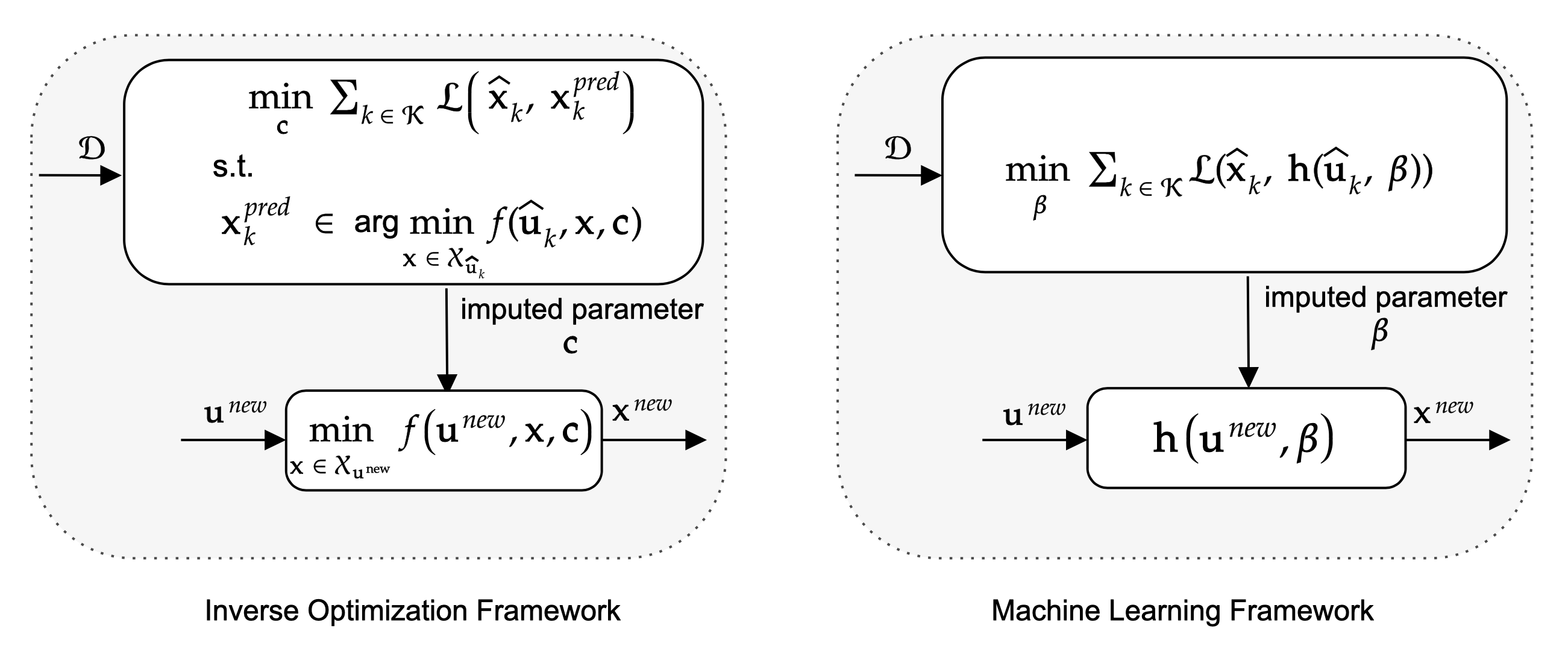}
\caption{IO and ML perspectives applied to an IO problem.}
\label{IO vs ML}
\end{figure}

Learning is the capability to acquire knowledge by extracting patterns from raw data \citep{Goodfellow-et-al-2016}; in this sense, IO fits within the broader conceptual framework of learning, i.e., we are acquiring knowledge about an optimization process that generated the data we observe. Applied to the IO problem described in Section \ref{sec:problem_def}, supervised ML would aim to find a vector-valued function $\mathbf{h}(\textbf{u}, \boldsymbol{\beta})$ from  $\textbf{u}$ (input variable or feature) to $\textbf{x}$ (target variable) such that the loss function $\pazocal{L}$ is  minimized \citep{russell2016artificial}. Once $\mathbf{h}(\textbf{u},\boldsymbol{\beta})$ is imputed, it can be used for prediction given new feature observations $\textbf{u}^{new}$. Figure \ref{IO vs ML} illustrates this observation: in the left panel, we shown that IO takes pairs of $(\hat{\bu}, \hat{\bx})$ (i.e., data) as input, imputes the value of $\bc$ (i.e., knowledge) and then, given $\bu^{new}$ and an imputed $\bc$, makes a prediction $\bx^{new}$; in the right panel, we show that the same happens in a typical supervised ML framework, with the difference being the actual models and methods used to impute the unknown parameters and, consequently, the learned model that is used for prediction.

The fundamental assumption in IO is that the observations are near-optimal or optimal solutions to an optimization problem, whereas general ML methods do not require such an assumption. Based on the description above, an IO problem can be seen as part of the same conceptual framework as machine learning -- the main difference being that in IO data is assumed to come from an optimization process. Viewed another way, we can say that IO assumes a strong prior that the input data is coming from an optimization process, and IO methods leverage this prior knowledge.  

\begin{figure}[htb]
\centering
\includegraphics[scale=0.14]{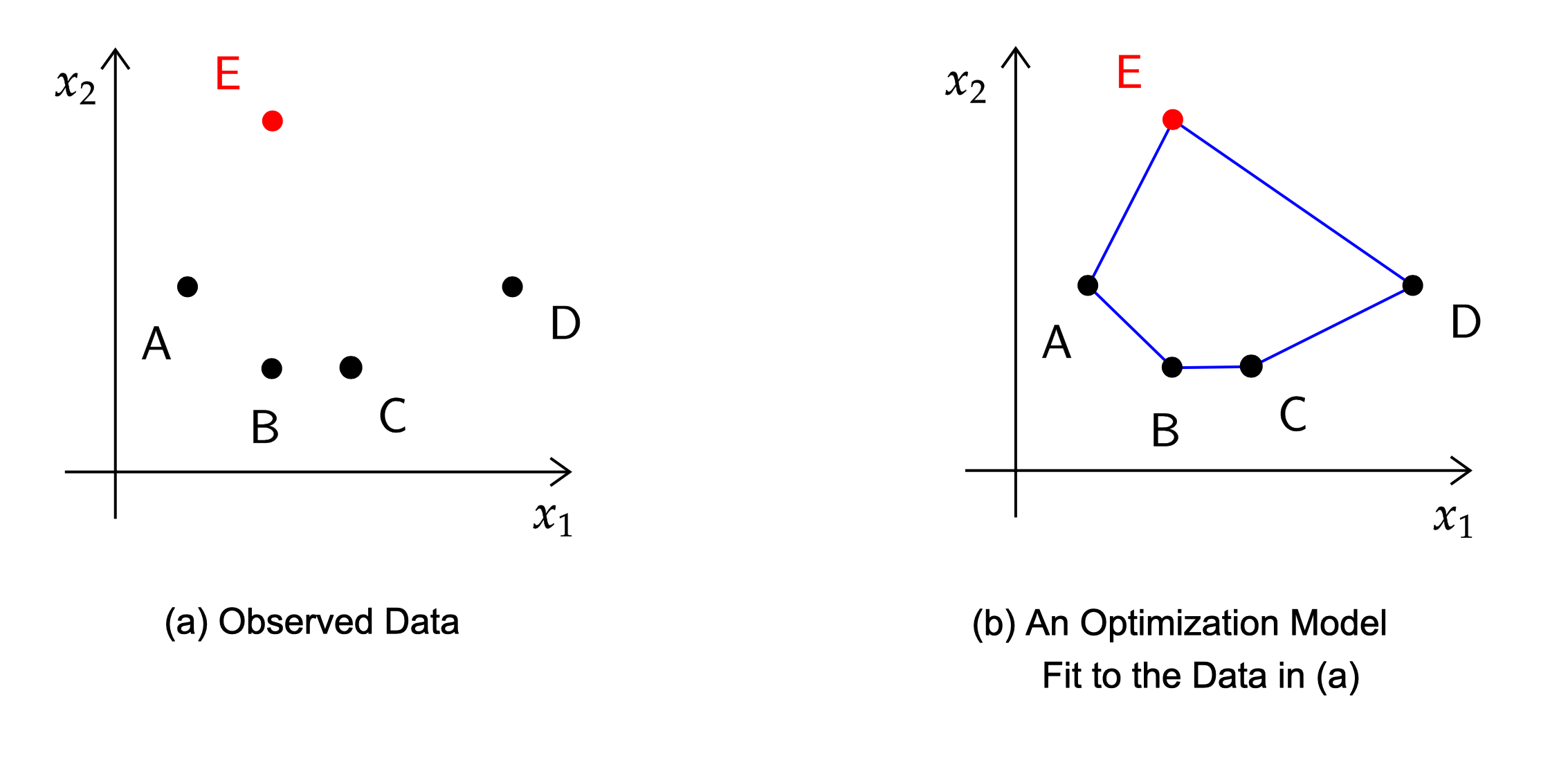}
      \caption{A data set and the fitted optimization model (points A to D are the training set, point E is the test set).}
    \label{fig:optimization_data}
\end{figure}

Consider Figure \ref{fig:optimization_data} (a): points A to D are given and our goal is to predict the next observation. Point $E$ is an unseen test point for which we want to evaluate the quality of prediction. Given no additional information, a natural idea would be to fit a curve using simple linear or polynomial regression. The training error for these models may be low but the test error would be high, which would normally signal over-fitting to the training data. However, here the problem would not be over-fitting --- rather, the issue is that the data is \emph{not} coming from a linear or a polynomial relationship between $x_1$ and $x_2$ and also not from a linear or polynomial relationship between a feature $u$ and each of the $x_i$. In reality, $A$, $B$, $C$ and $D$ are optimal solutions to the parametric optimization problem $P(u):\max \; (-6+5u)x_{1}-(3+10u)x_{2}$ subject to the convex hull of $(A, B, C, D, E)$, shown in Figure \ref{fig:optimization_data} (d), with $P(-0.96 \leq u \leq 0.16)=A$, $P(-0.16 \leq u\leq 1.25)=B$, $P(1.25 \leq u \leq 17.7)=C$, $P(17.7 \leq u \leq 20)=D$ and $P(-20 \leq u \leq -0.96)=E$. This example was generated using the Wolfram Parametric Linear Programming app \citep{Mathematica}. 

Just like we aim to fit a linear regression model when we conjecture the relationship is linear, or a polynomial one when we conjecture it is polynomial, we can fit an optimization model when we think the data is coming from an optimization problem. Prior knowledge of the process that generated the data can have a substantial impact on the ability to fit a good model; the above example illustrates specifically the effect of incorporating an optimization model prior. Next, we compare the predictive performance of IO methods, which assume an optimization prior, and ML methods, which do not, when data is generated from parametric convex optimization problems. 

\section{Experiments}\label{sec:experimentalSetup}

\subsection{Preliminaries}
Our goal is to evaluate `out-of-the-box' approaches that can be easily implemented in practice: a KKT-based IO model \citep{Keshavarz11}, random forest \citep{breiman2001random}, support vector regression \citep{vapnik95} and Gaussian process regression \citep{rasmussen2004gaussian}. We use the scikit-learn \citep{scikit} implementation of the ML methods. \citet{Keshavarz11}'s KKT-based IO model is solved using \emph{Concert Technology} of IBM ILOG CPLEX v.12.7.0. All experiments are run on an HP server with 20 Intel(R) Xeon(R) CPU E5-2687W v3 @ 3.10GHz processors and 512 GB RAM under Linux environment.
 
To show the value of prior knowledge for IO, we consider two types of IO models. In \emph{IO-perfect}, we have perfect information about the \emph{form} of the unknown objective function (but do not know its parameters). In \emph{IO-imperfect}, the objective function is misspecified, and so differs from the true one due to modeling error. 


\noindent
\begin{algorithm}[H]
\caption{Data Generation}\label{euclid}
\begin{algorithmic}
\For  {$k \gets 1$ to $K$}
\State generate $\textbf{u}_{k} \textrm{ from }  \textrm{uniform dist'n}$
\State solve $\mathbf{FOP}_{true}(\bu_k, \bc^{true})$ and get ${\textbf{x}}_{k}^{true}$ 
\State $k \gets k+1$
\EndFor
\State $\pazocal{D} \gets \{(\hat{\textbf{x}}_{1},\hat{\textbf{u}}_{1}), (\hat{\textbf{x}}_{2},\hat{\textbf{u}}_{2}), ..., (\hat{\textbf{x}}_{K},\hat{\textbf{u}}_{K})  \}  $ \vspace{0.2cm}
\end{algorithmic}
\label{algo}
\end{algorithm}

Pairs $(\hat{\textbf{u}}_k, \hat{\textbf{x}}_k)$ are generated using Algorithm \ref{algo} and divided into a training set and a test set. We use leave-one-out cross-validation on the training set for the ML methods: we pick one point as our validation set, fit a model to the remaining data and evaluate the error on the single held-out point. A validation error estimate is obtained by repeating this procedure for each of the training points and averaging the results. Doing so with several hyperparameter settings enables us to choose the hyperparameter setting with the lowest average cross-validation error. All of the ML experimental results presented below show the error on the test set for the best hyper-parameter setting found during the leave-one-out cross-validation procedure.

We measure the deviation between the predicted $\textbf{x}$ and $\textbf{x}^{true}$ in the test set (of size $M$) using mean relative error (MRE):
\begin{equation}
    MRE = \frac{1}{M} \sum _{m=1}^{M} \frac{ \| \textbf{x}- \textbf{x}^{true} \|_2 }{ \| \textbf{x}^{true} \|_2}.
\end{equation}
 
\subsection{Experiment 1: Utility Function Estimation}
Consider imputing a customer's unknown utility function given observations of their past purchases \citep{Keshavarz11, Barmann17}. We assume that customers solve an internal optimization problem in the course of purchasing goods to maximize their utility function:
\begin{equation}\label{eq:UFOP}
\begin{alignedat}{1}
UFOP (\textbf{\textbf{p}}) : \underset{\textbf{x}}{\text{minimize}} & \quad \textbf{p}^{T}\textbf{x} -U(\textbf{x})  \\
\text{subject to}
 & \quad \bx \ge \bzero,
\end{alignedat}
\end{equation}
\noindent
where $\textbf{x}\in \mathbb{R}^{n}$ and $\textbf{p}\in \mathbb{R}^{n}$ are the vectors of consumer demand and the associated price, respectively. Similarly to \citet{Keshavarz11}, we assume that the function $U(\textbf{x})$ is a concave utility function which is non-decreasing over $[0, x^{\max}]$, $x^{\max}$ being the maximum demand level found from past observations. Since $U(\bx)$ is concave, the objective function of (\ref{eq:UFOP}) is convex. 

We generate $p_i \sim U[5,20]$. The true utility function is $U_{true}(\textbf{x})= \sum_{i=1}^{n} \sqrt{x_{i}}$. We use $U_{perfect}(\bx)= c_{i}\sqrt{x_{i}}$ in implementing \emph{IO-perfect}  (i.e., the model is correctly specified so that all observations can be made optimal solutions for their corresponding $p_{i}$). For \emph{IO-imperfect}, we use $U_{imperfect}(\textbf{x})= \sum_{i=1}^{n} q {x}_{i}^{2}+ 2r{x_i}$.

Based on the results of our hyper-parameter search during leave-one-out cross-validation, we use the scikit-learn implementation of random forest with \texttt{n\_estimators=50}, \texttt{max\_depth=None}, \texttt{max\_feature=None}, support vector regression with  \texttt{c=0.1}, \texttt{gamma="auto"}, \texttt{kernel=RBF} and Gaussian process with \texttt{kernel=RBF}, \texttt{alpha=1e-10}, \texttt{normalized\_y=False}, \texttt{optimizer="fmin\_1\_bfgs\_b"} and \\ \texttt{n\_restarts\_optimizer=10}. 

\begin{figure}[H]
\centering\includegraphics[width=0.6\linewidth]{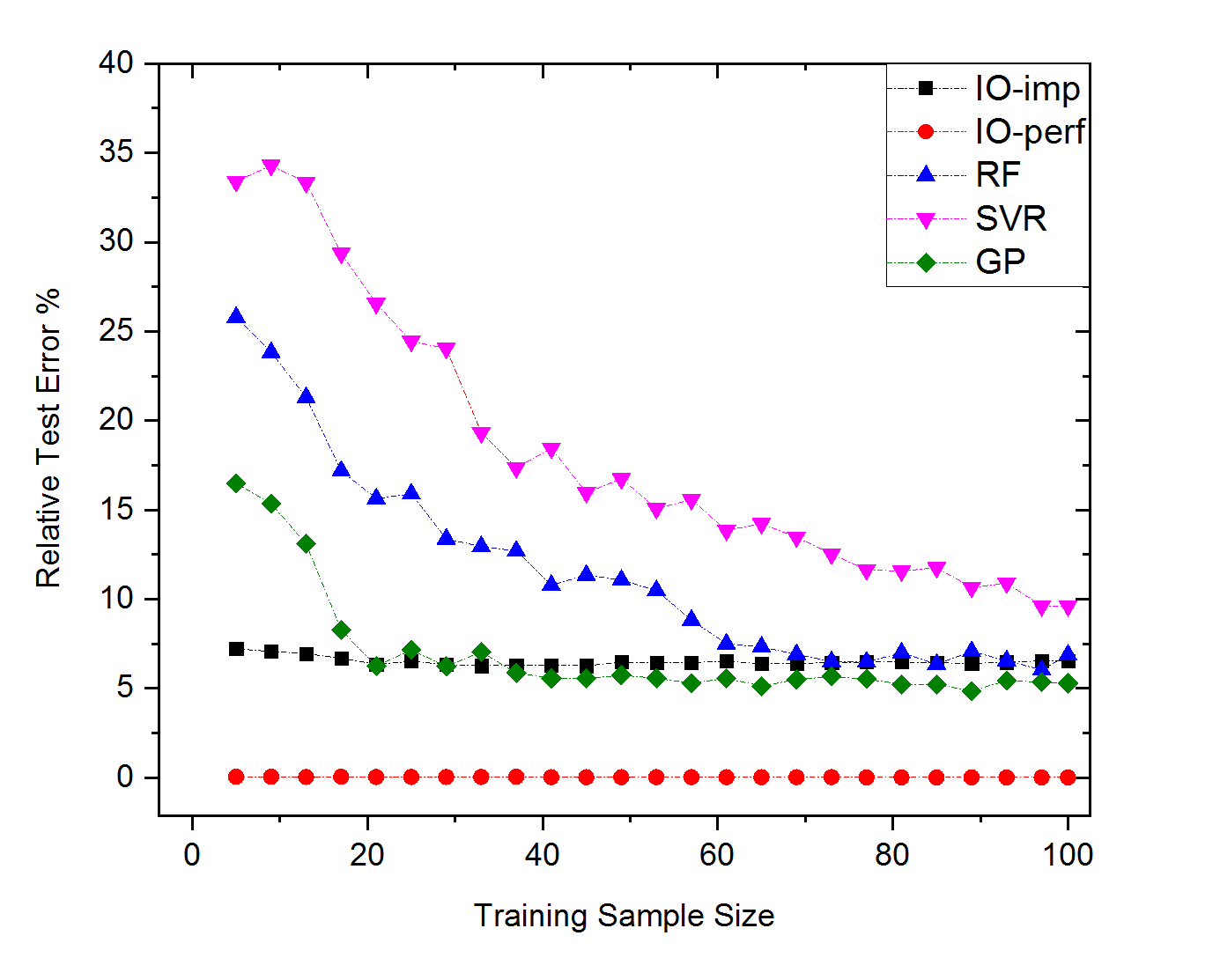}
\caption{Comparing \emph{IO-imperfect}, \emph{IO-perfect}, RF, SVR, and GP for utility function estimation; each point is an average over 15 randomly-generated problems, $p_i \sim U[5, 20]$.}
\label{keshavarzfig}
\end{figure}

Figure \ref{keshavarzfig} shows how the relative test error of \emph{IO-perfect}, \emph{IO-imperfect}, RF, SVR and GP changes when we increase the training size. Whereas RF, SVR and GP performance improves, the performance of \emph{IO-imperfect} and \emph{IO-perfect} is not affected. GP achieves performance similar to that of \emph{IO-imperfect} with only 20 observations, while RF and SVR get close to \emph{IO-imperfect} at 60 and 100 points, respectively. The ML methods perform well because the underlying optimization process (resulting from model (\ref{eq:UFOP})) is constrained by only non-negativity constraints and is parametric only in the objective function -- the relationship between the input and output of this process can be captured (learned) without knowledge of the underlying optimization model. In the next experiment, we aim to identify IO problems which would be more difficult for these ML methods and would require knowledge of the underlying optimization structure to make good predictions. 

\subsection{Experiment 2: Randomly-Generated Parametric Optimization Problems}
Using the POP toolbox \citep{POP}, we generate POPs of the following form:
\begin{equation}\label{eq:TPOP}
\begin{alignedat}{1}
TPOP (\textbf{u}) : \underset{\bx}{\text{min}} & \quad (\mathbf{Q}\bx+\mathbf{H}\bu+\bc)^{T}\bx\\
\text{s.t.}
 & \quad \bA\bx \leq \bb+ \mathbf{F}\bu\\
 & \quad \bx \in \mathbb{R}^{n} \\
 & \quad \bu \in \mathbb{U}:= \{\bu \in \mathbb{R}^q \hspace{0.1cm} | \hspace{0.1cm} \mathbf{CR}_{A}\bu \leq \mathbf{CR}_{b} \}, \\
\end{alignedat}
\end{equation}

\noindent
where $\mathbf{Q} \in \mathbb{R}^{n\times n} \succ 0$, $\mathbf{H} \in \mathbb{R}^{n\times q}$ and $\bc \in \mathbb{R}^{n\times 1}$. In the inequality constraints, $\bA \in \mathbb{R}^{m \times n}$, $\bb \in \mathbb{R}^{m \times 1}$ and $\mathbf{F} \in \mathbb{R}^{m \times q}$. 
The last constraint of (\ref{eq:TPOP}) restricts the parameter space, i.e., the admissible values of external parameter (feature) $\bu$, where $\mathbf{CR}_{A} \in \mathbb{R}^{2q \times q}$ and $\mathbf{CR}_{b} \in \mathbb{R}^{2q \times 1}$ allow for representations of constraints on $\bu$ (in particular, lower and upper bounds -- hence the number of rows is $2q$), and 
$q$ denotes the dimensionality of $\bu$. We assume $q = 2$ in our experiments.

\begin{figure}[H]
    \centering 
   \includegraphics[width=1.00\linewidth]{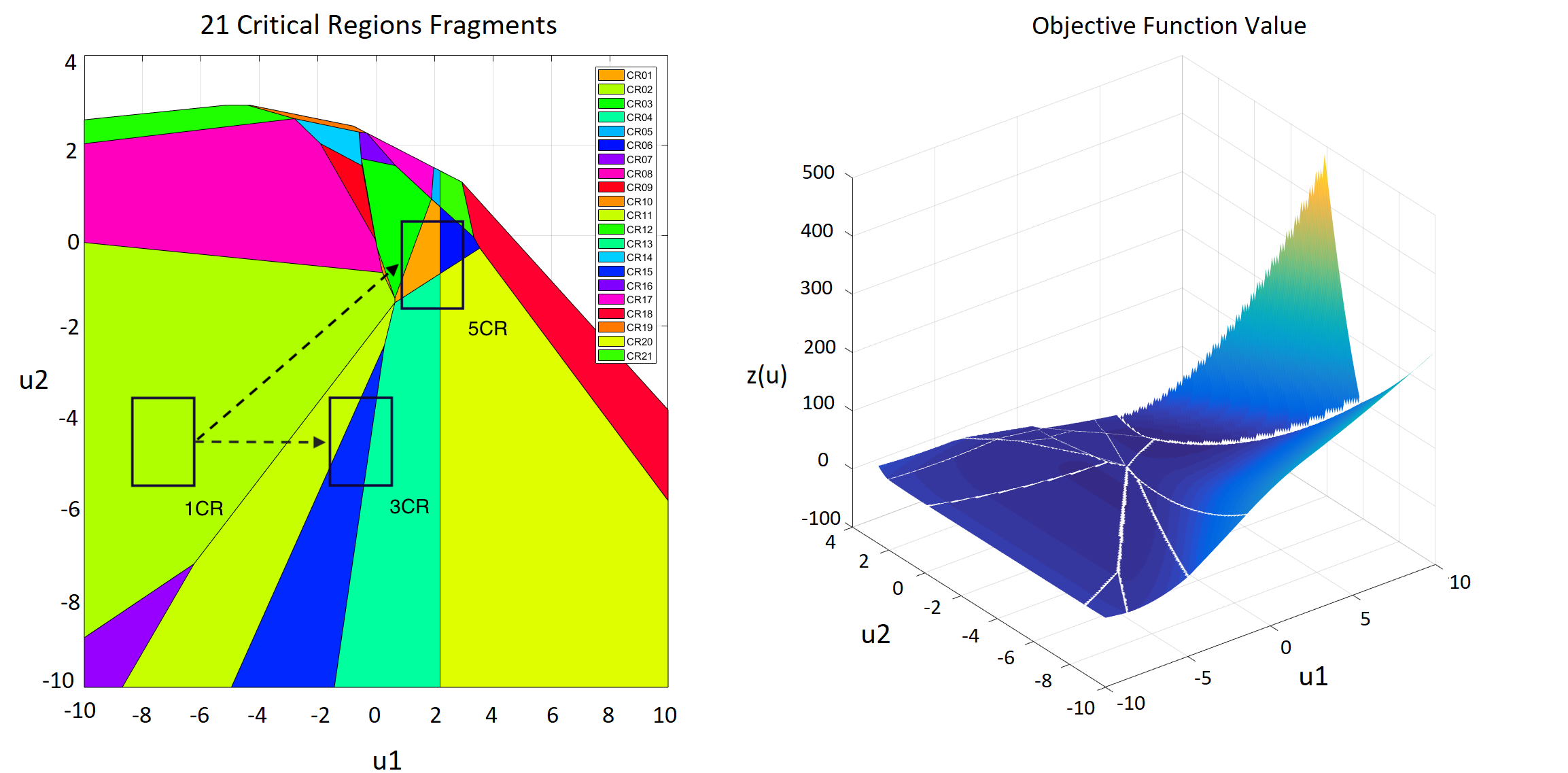}
    \caption{A POP instance defined over two parameters $u_1$ and $u_2$ generated using the POP toolbox of \citet{POP}. The left panel shows 21 critical regions of the parameter space, while the right panel shows the piecewise quadratic optimal objective function $z(\bu)$. Each critical region on the left corresponds to a quadratic `piece' of $z(\bu)$ on the right.}
    \label{CR distribution}
\end{figure}

To study (forward) POPs, it is standard to view the space of feasible parameters $\bu$ as being partitioned into polyhedrons called \emph{critical regions}, each of which is uniquely defined by a set of optimal active constraints \citep{Ahmadi18,POP}. It is known that if $\mathbf{Q}$ is positive definite and $\mathbb{U}$ is convex, the optimal objective function $z(\bu)$ is continuous and piecewise quadratic, while if the problem is reduced to a multi-parametric linear program, then the optimal objective function $z(\bu)$ is continuous, convex, and piecewise affine (see \citet{POP} and the references therein). Figure \ref{CR distribution} shows an example of one of the instances we generate using the POP toolbox of \citet{POP} with 21 critical regions and a two-dimensional $\bu$; the critical regions are shown in the left panel while the corresponding optimal value function $z(\bu)$ is shown on the right. 

To be able to control the number of critical regions in our experiments, we manually pre-process the instances generated from the toolbox to define problems over smaller regions of the parameter space with a particular number of critical regions. For instance, from the example in Figure \ref{CR distribution}, we can generate a one-region, a three-region and a five-region problems by selecting the regions indicated by the black rectangle borders. Throughout this process, care is taken to select intervals where each critical region has close to a ``fair share'' of the region; in preliminary experiments, it was noted that the difficulty of learning over a parameter space interval that covers, say, five regions but is dominated by one particular region is similar to the case of one critical region (more discussion on the effect of critical regions on learning performance appears later in this section). 

\paragraph*{Varying the Training Set Size (Data Efficiency)}
We compare the predictive performance of \emph{IO-perfect}, \emph{IO-imperfect}, SVR, RF and GP varying the training size from 1 to 500 on 20 randomly generated POPs of type \eqref{eq:TPOP} with two parameters $u_1$ and $u_2$, with 1 and 3 critical regions. For \emph{IO-imperfect} we use the mean of selected intervals of $\bu$ in the $(\textbf{H} \bu)^{T} \bx$ term in the objective function instead of a parametric term.

\begin{figure}[htb]
    \hspace{-3.3em}\subfloat[One Critical Region]{{\includegraphics[trim=0cm 0cm 0cm 2cm, clip, width=0.59\linewidth]{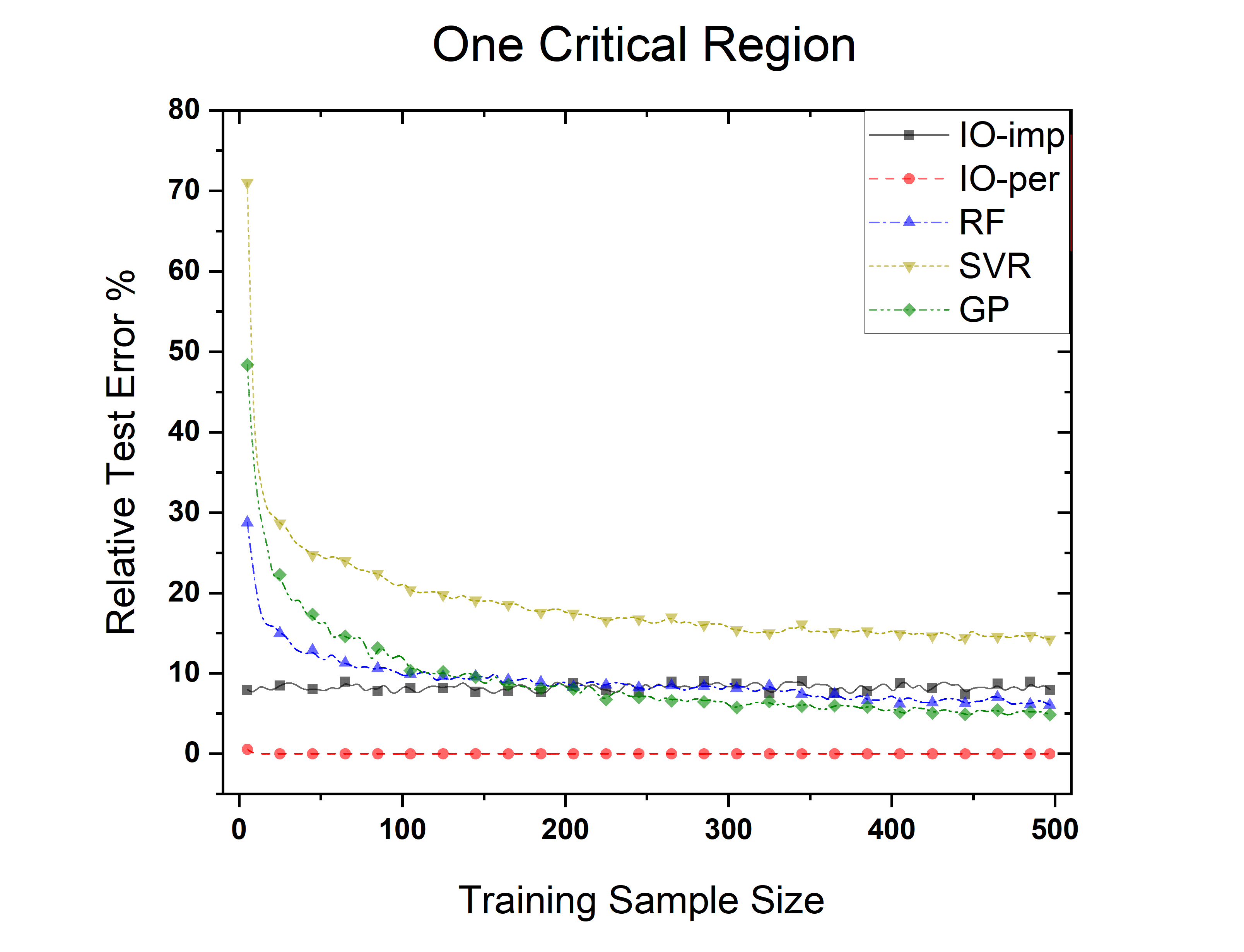} }}%
    \subfloat[Three Critical Regions]{{\includegraphics[trim=0cm 0cm 0cm 2cm, clip, width=0.59\linewidth]{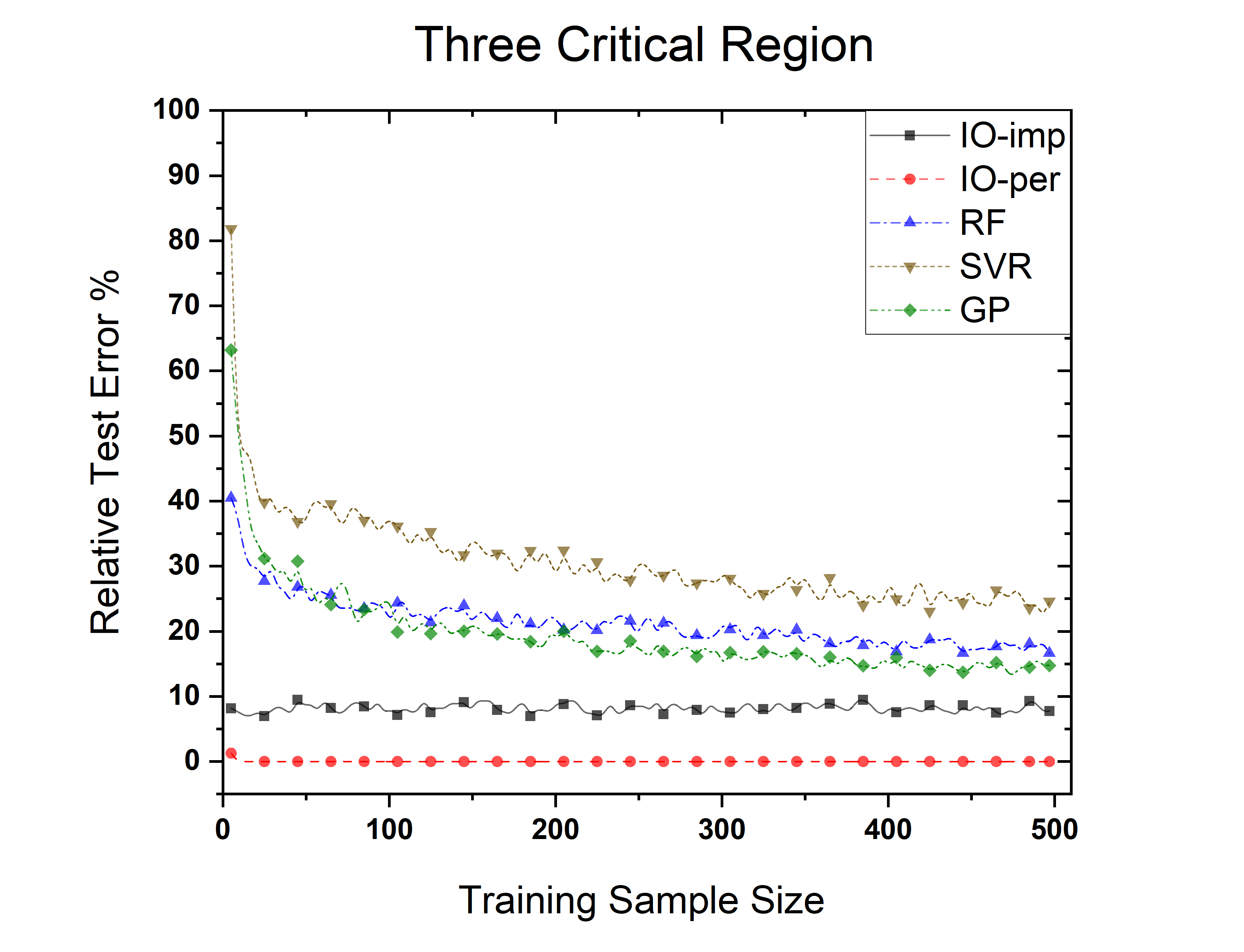} }}%
      \caption{Comparing \emph{IO-perfect}, \emph{IO-imperfect}, SVR, RF and GP varying training sample size and critical regions.}%
      \label{CR1}
\end{figure}

Figures \ref{CR1}(a) and \ref{CR1}(b) illustrate the predictive (i.e., test set) performance of these methods (averaged over 20 instances) given data sets from one critical region and three critical regions, respectively. As expected, ML performance improves with increasing training set size in both cases. On the contrary, IO methods are not affected: \emph{IO-perfect} has nearly zero error while \emph{IO-imperfect} incurs around 8\% error, regardless of the training set size. Since IO is able to achieve low prediction error with a much smaller training size, we see that IO methods with correctly specified priors can be more data efficient than ML algorithms for problems where data is generated by an optimization process.

Because the IO approach we chose solves an optimization model, no matter how many points are given, \emph{IO-perfect} is able to recover the parameters that make the observations optimal while \emph{IO-imperfect} will maintain the same level of mis-specification. Insensitivity to the training set size is consistent with \citet{aswani2018inverse}'s results, which showed that \citet{Keshavarz11}'s method is in general not risk consistent, i.e., it is not guaranteed to asymptotically provide the best possible predictions. However, we emphasize that despite being risk inconsistent, given the right prior, the \citet{Keshavarz11} approach performs well, and can be the method of choice due to its data efficiency. 

The test error for ML methods increases for all training set sizes as we move from one critical region (Figure \ref{CR1}(a)) to three critical regions (Figure \ref{CR1}(b)). Comparing Figures \ref{CR1}(a) and \ref{CR1}(b), the curves for GP and RF in \ref{CR1}(a) intersect \emph{IO-imperfect} after 200 points, but in \ref{CR1}(b) they require more points to be comparable to \emph{IO-imperfect}. This observation suggests that the number of critical regions in the parameter space may be a factor that makes these problems more difficult for standard ML approaches, motivating the next experimental comparison. 

\begin{figure}[H]
    \centering
   \includegraphics[width=1.0\linewidth]{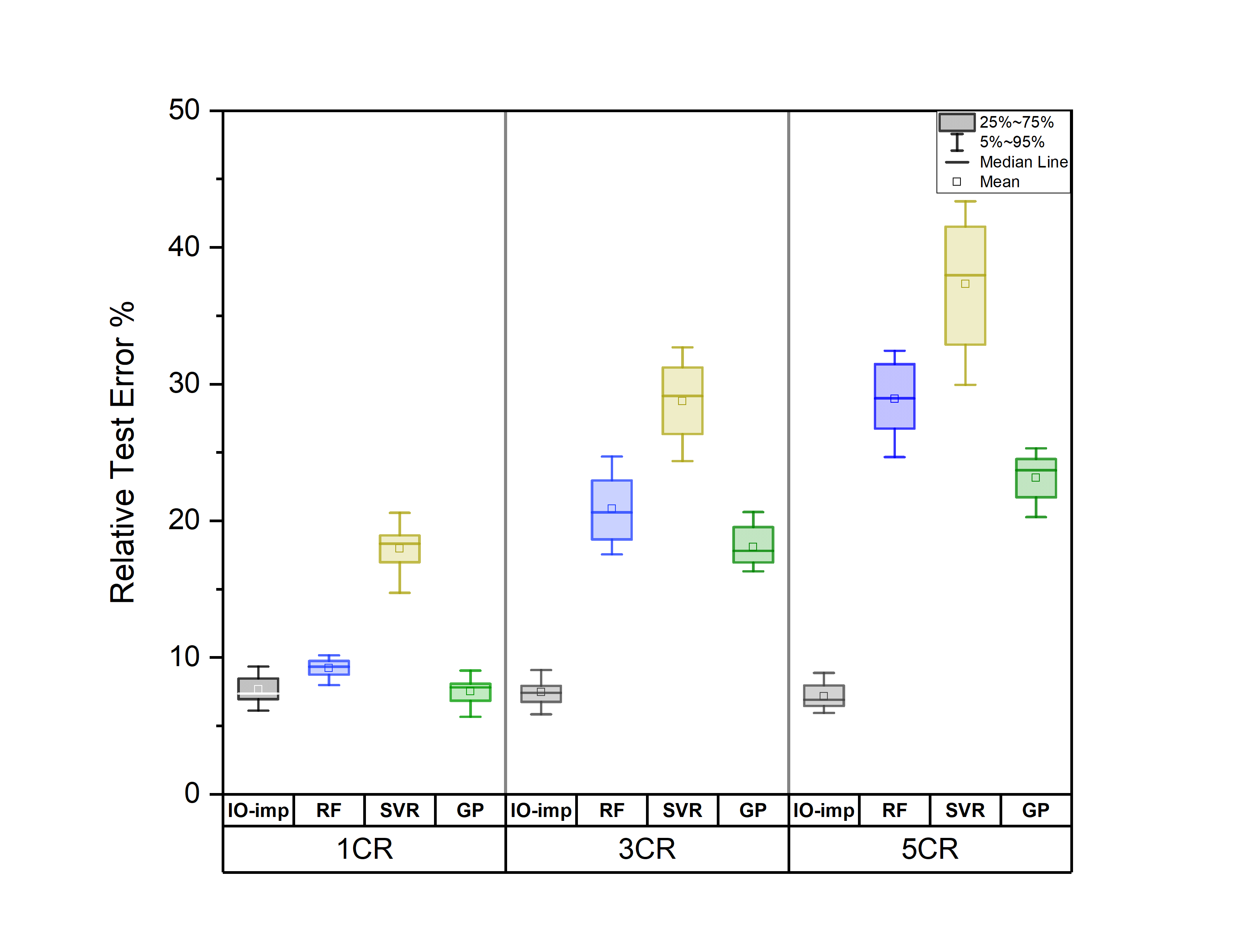}
    \caption{Comparing \emph{IO-imperfect}, SVR, RF and GP as the number of critical regions increases; training set size = 200, test set size = 200.}
    \label{error vs CRs}
\end{figure}

\paragraph*{Increasing the Number of Critical Regions} 
Figure \ref{error vs CRs} shows the average predictive performance of \emph{IO-imperfect} and the three ML methods over 20 randomly generated POPs given a training set of 200 observations. For each of the 20 problem instances, we choose one interval of $\bu$ with one critical region, one with three critical regions and one with five, keeping the chosen range of both $u_1$ and $u_2$ equal to two, as in Figure \ref{CR distribution}. We can see that \emph{IO-imperfect} is not sensitive to the number of critical regions while for all employed ML methods the test error increases with the number of critical regions. As mentioned above, every critical region is a polyhedron of $\bu$ values over which the optimal objective function $z(\bu)$ is quadratic (a quadratic `piece' of the overall piecewise quadratic function) while the optimizer function $\bx^*(\bu)$ is piecewise-affine \citep{POP}. It is not surprising that some of the ML methods can perform well with 200 points over one critical region as they only have to learn one quadratic function. When the number of critical regions increases, the task for the chosen ML methods becomes more challenging: they have to learn multiple quadratic functions without any prior knowledge about when the transition from one region to another will occur. We do not rule out the possibility of creating new, custom learning algorithms to overcome this issue but doing so is beyond the scope of this paper. 

\begin{figure}[H]
  \centering
    \includegraphics[width=1.00\linewidth]{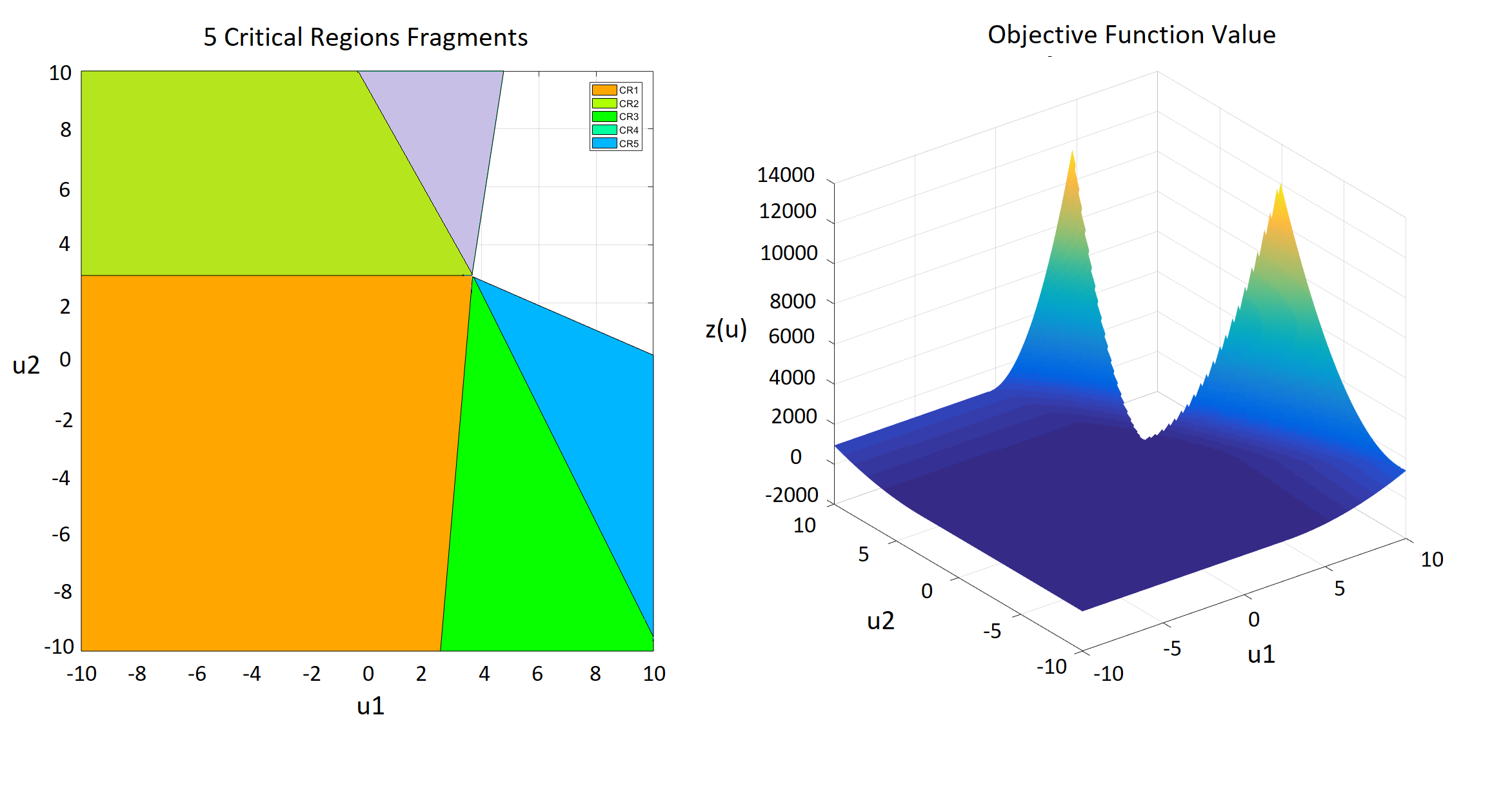}
    \caption{Problem \eqref{eq:FOP2} critical regions (left) and  parametric objective function (right) generated using the POP Toolbox of \citet{POP}.}
    \label{5 CR instance}
\end{figure}

\paragraph*{Varying the Functional Dependence of the Objective Function on $\textbf{u}$} Next, we investigate how the relationship between the TPOP objective function and the feature $\textbf{u}$ affects the performance of ML and IO. We investigate this question in the context of one randomly generated POP (see Figure \ref{5 CR instance}): 

\begin{equation}\label{eq:FOP2}
\begin{alignedat}{1}
z(\textbf{u}) : \underset{\bx}{\text{minimize}} & \quad 1.3040 x_{1}^{2}+(1+u_{1})x_{1}+ 19.4545x_{2}^{2}+(u_{2}-u_{1}+1)x_{2}  \\
\text{s.t.}
 & \quad 0.2294 x_{1} \leq 2.5237- 0.9733 u_{2}\\
 & \quad 0.1890 x_{2} \leq 4.2679 - 0.8658 u_{1}- 0.4634u_{2}\\
 & \quad 0.5436 x_{1} - 0.5889 x_{2} \leq 2.8088 - 0.4757 u_{1} - 0.3624 u_{2} \\
 & \quad 0.2210 x_{1} \leq 3.2535 - 0.9753 u_{1}.
\end{alignedat}
\end{equation}

We consider three objective functions for problem \eqref{eq:FOP2} defined through functions $\Psi_1(\textbf{u})$ and $\Psi_2(\textbf{u})$, which set the coefficients of $x_1$ and $x_2$: 
\begin{equation}
   f= 1.3040x_{1}^{2}+\Psi_{1}(\textbf{u}) x_{1} + 19.4545x_{2}^{2}+ \Psi_{2}(\textbf{u}) x_{2}.
\end{equation}
\noindent 
We start with a linear $\Psi(\bu)$ followed by changing it to a multi-variate rational function and subsequently increasing the degree of the polynomials forming the rational function to create a relationship that could lead to a more challenging learning problem. 

As shown in Table \ref{obju}, the test error of GP, SVR and RF increases as we move from the `simple' (linear) dependence on $\bu$ to making the coefficient of $\bx_1$ and $\bx_2$ multi-variate rational functions of $\bu$. Since ML methods aim to find a mapping from $\textbf{u}$ to $\textbf{x}$, a more complex $\Psi(\textbf{u})$ indeed makes learning more difficult. On the other hand, IO is not sensitive to the nature of $\Psi(\textbf{u})$ as this function is just the coefficient of $\textbf{x}$ in the IO mathematical model. 

\begin{table}[H] \caption{Comparing the predictive performance of ML and IO varying the form of the true objective function, with $\hat{u}_{1}\sim U(4,6)$ and $\hat{u}_{2}\sim U(-6,-4)$. Recall that $K$ is training set size. The test set size is set to be equal to the training set size.}
\centering
\begin{tabular}{llcc} 
\hline
True Objective Function & Method
&\makecell{MRE \% \\ ($K=100$)} &\makecell{MRE \% \\ ($K=200$)}      \\
 \hline \\
\multirow{5}{*}{\makecell{$1.3040x_{1}^{2}+(1+u_{1})x_{1}$\\ $+19.4545x_{2}^{2}+(u_{2}-u_{1}+1)x_{2}$}} & \emph{IO-perfect} & $\sim 0$ &$\sim 0$  \\
                  & \emph{IO-imperfect}  &3.85& 3.79    \\
                  & GP & 9.75 & 6.03  \\
                  & RF & 6.75 &5.74\\
                  & SVR & 19.03 & 18.04\\ \\ \hline \\
\multirow{5}{*}{\makecell{$1.3040x_{1}^{2}+ \frac{1+u_{1}}{1-u_{1}}x_{1}$ \\ $+19.4545x_{2}^{2}+ \frac{(u_{2}-u_{1}+1)^{3}}{(u_{2}-1)^{2}}x_{2}$}} & \emph{IO-perfect} & $\sim 0$  & $\sim 0$ \\
 
                  & \emph{IO-imperfect}  &7.80&   7.11  \\
                 & GP &10.85 & 8.59 \\
                & RF &9.77 & 7.07\\
                  & SVR &24.94 & 23.30\\ \\  \hline \\

\multirow{5}{*}{\makecell{$1.3040x_{1}^{2}+ \frac{u_{1}}{(1-2u_{1})^{4}}x_{1}$ \\ $+19.4545x_{2}^{2}+ \frac{(u_{2}-u_{1})^{3}}{(3u_{2}-5)^{5}}x_{2}$}} & \emph{IO-perfect} & $\sim 0$ &$\sim 0$    \\
 
                  & \emph{IO-imperfect}  &7.78 &   7.62  \\
                 & GP &15.28& 10.29 \\
                & RF &22.24 & 17.75 \\
                  & SVR &38.76 & 34.80 \\ \\  \hline
\end{tabular}
\label{obju}
\end{table}

\paragraph*{Investigating the Impact of Correctness of Prior Knowledge}
Above, we show that \emph{IO-perfect} and \emph{IO-imperfect} are insensitive to problem characteristics which significantly impact the performance of ML (training set size, number of critical regions and the nature of the dependence on $\bu$). However, a difference between \emph{IO-perfect} and \emph{IO-imperfect} was observable. To investigate the impact of the correctness of objective function prior, we consider five objective function choices for problem (\ref{eq:FOP2}) which encode different degrees of correctness (or mis-specification) of a prior. 

In Table \ref{value}, given $\pazocal{D}$, we aim to impute $c_{1}$ and $c_{2}$. The first objective function is the exact objective function of problem \eqref{eq:FOP2}, i.e., the \emph{IO-perfect} method. The subsequent functions deviate from the true function in the coefficients of $x_1$ and $x_2$ (but the goal is still to impute $c_1$ and $c_2$ from the data generated by the true objective function of problem \eqref{eq:FOP2}). Besides parametric objective functions (i.e., objective functions defined over $\textbf{u}$), we also consider two non-parametric choices. In the first non-parametric objective, we use the mean of $1+\hat{u}_{1}$ and $\hat{u}_{2}-\hat{u}_{1}+1$ from the given training set $\hat{u}_1$ and $\hat{u}_2$; in the second non-parametric function, we completely eliminate the linear terms. 
 
We see that even slight differences in the function increase IO test error, although for the parametric functions the difference is small. The non-parametric mean-based function also does not lead to a significant increase, but this could change if a larger range of $\hat{u}$ is considered. Removing any encoding of the parametric parts of the objective function drastically increases the prediction error from 3.88\% to 52.28\%. These observations confirm the sensitivity of IO to the correctness of objective function prior. 

\begin{table} 
\centering
\caption{Different parametrizations of the objective function of problem \eqref{eq:FOP2} and the corresponding test error, given $\hat{u}_{1}\sim U(4,6)$ and $\hat{u}_{2}\sim U(-6,-4)$.}
\begin{tabular}{ccc}
\hline
                &  \makecell{Objective Function} & \makecell{MRE \% \\ ($K=100$)}      \\
 \hline
\multirow{3}{*}{Parametric} & $c_{1}x_{1}^{2}+(1+u_{1})x_{1}+c_{2}x_{2}+(u_{2}-u_{1}+1)x_{2}$  &$\sim$ 0    \\
                  &$c_{1}x_{1}^{2}+x_{1}+c_{2}x_{2}+(u_{2}-u_{1})x_{2}$  & 0.38    \\
                  &$c_{1}x_{1}^{2}+x_{1}+c_{2}x_{2}-u_{1}x_{2}$ & 3.56 \\ 
\multirow{3}{*}{Non-Parametric} &  &  \\
                  &$c_{1}x_{1}^{2}+ 6 x_{1}+c_{2}x_{2}-9 x_{2}$  & 3.88  \\
                  & $c_{1}x_{1}^{2}+c_{2}x_{2}$ & 52.28  \\ \hline
\end{tabular}

\label{value}
\end{table}

\subsection{Insights}
Our experiments with randomly-generated POPs demonstrate the importance of four characteristics as determinants of the difficulty of objective function learning: (i) size of the training set, (ii) nature of the dependence of the optimization problem on the external parameters, (iii) level of confidence with regards to the correctness of the optimization prior, (iv) number of critical regions in the parameter space of a POP. 
    
Increasing the size of the training set favours ML methods but not necessarily the IO methods. The contrast in performance between the problem in experiment 1 and the randomly-generated POPs suggests that ML methods are more likely to be successful on problems with no or few (parametric) constraints. Increasing the complexity of the relationship between $\bu$ and the objective function also makes the problem more difficult for ML. On the other hand, the performance of IO is strongly dependent on the correctness of the objective function prior. 
    
Perhaps the most important characteristic determining the difficulty of the objective function learning tasks is the number of critical regions. Over 20 randomly-generated POPs, we found a substantial increase in the relative test error for all three ML methods as the number of critical regions increased. In fact, we conjecture that the underlying reason why adding parametric constraints and/or increasing the complexity of the dependence on $\bu$ made learning more difficult is because doing so also increased the number of critical regions. 

Given these observations, a summary of recommendations of when to use ML or IO for learning the coefficients of a convex POP is needed. We summarize our recommendations in Table \ref{tab:recommend}. In each combination of a criterion and its value (low/high or small/large), we state whether we expect the classic ML methods or an IO method to be a better choice (or at least a better starting point) for solving the learning problem, while holding other criteria fixed. For example, when training set size is small, we expect IO to be a better option due to its data efficiency. 

\begin{table}[h] 
\centering
\caption{Problem Characteristics and Suggested Methods.}
\begin{tabular}{ccc}
\hline
                &  \makecell{Low/Small} & \makecell{High/Large}      \\
 \hline
Size of Training Set & IO & ML    \\
Dependence on External Parameter & ML & IO \\
Confidence in Correctness of Prior & ML & IO \\
Number of Critical Regions & ML & IO \\ \hline
\end{tabular}
\label{tab:recommend}
\end{table}
    
The relative performance of the methods considered may not be the same for other types of learning problems where data was generated by an optimization problem (e.g., learning constraints, learning in discrete optimization problems). However, we conjecture that the analysis of the underlying structure of the value function or the optimizer function in terms of the external parameters $\bu$ (i.e., analysis of the problem in terms the critical regions) will be important in both gaining additional understanding of the challenges of learning from optimization data and developing more sophisticated learning methods for such problems.

\section{Conclusion} \label{conclusion}
In this paper, we view inverse optimization as a problem of learning from decisions that are made through an unknown optimization process. We specifically focus on the problem of learning a convex objective function of a parametric optimization problem. We experimentally compare the predictive performance of an inverse optimization method with perfect and imperfect priors with three well-known machine learning algorithms: support vector regression, random forest and Gaussian processes. While we show that some inverse optimization problems can be tackled through classic machine learning approaches, we highlight the need for sophisticated inverse optimization models for problems where at least one of the following characteristics holds: (i) the size of the training set is small, (ii) both the constraints and the objective function of the problem in question are dependent on an external parameter (feature), particularly when that dependence is non-linear, (iii) we have high confidence that our knowledge of the parametric nature of the constraints and objective is correct, and (iv) the number of critical regions of the POP is large with no one region dominating. We believe that these observations provide practitioners with guidance on when to consider employing inverse optimization instead of, or in addition to, classical machine learning. 

\appendix
\section*{Appendix: KKT-based Inverse Optimization Model} The \citet{Keshavarz11} IO model for (\ref{eq:FOP1}) is 
\begin{equation}
\nonumber
\begin{alignedat}{1}
\mathbf{IO} (\pazocal{D}) : \underset{\blambda, \textbf{c}}{\text{minimize}} & \quad
\sum_{k \in \pazocal{K}} \phi (\textbf{r}_{k}^{stat},\textbf{r}_{k}^{comp})\\
\text{s.t.}
 & \quad {r^{stat}}_{k}^{j} =\frac{\partial f(\textbf{u}, \textbf{x}, \textbf{c})}{\partial{x_{j}} }+\sum_{i= 1}^{m}\lambda_{i}^{k} \frac{\partial g_{i}(\textbf{u}, \textbf{x})}{\partial {x_{j}} } \hspace{0.2cm} \quad \forall j\in \pazocal{J}, k \in \pazocal{K}\\
 & \quad {r^{comp}}^{i}_{k}  = -\lambda^{i}_{k}  g_{i}(\textbf{u}, \textbf{x}) \quad \quad \hspace{2.4cm} \forall i\in \pazocal{I},  k\in \pazocal{K}\\
& \quad \blambda \ge \textbf{0},
\end{alignedat}
\end{equation}
where $\phi$ is some norm, $\textbf{r}^{stat}$ and $\textbf{r}^{comp}$ represent the complementary slackness and stationarity residuals, respectively, of the KKT conditions; $\blambda$ is the dual variable. The missing objective function coefficient vector $\textbf{c}$ is a decision variable and pairs of $(\hat{\textbf{x}}_{k}, \hat{\textbf{u}}_{k})$ are the inputs.

\section*{Acknowledgement(s)}

The authors are supported by the Natural Sciences and Engineering Research Council
of Canada.

\bibliographystyle{apacite}
\bibliography{interactapasample}

\begin{thebibliography}{}

\bibitem [\protect \citeauthoryear {%
Ahmadi-Moshkenani%
, Johansen%
\BCBL {}\ \BBA {} Olaru%
}{%
Ahmadi-Moshkenani%
\ \protect \BOthers {.}}{%
{\protect \APACyear {2018}}%
}]{%
Ahmadi18}
\APACinsertmetastar {%
Ahmadi18}%
\begin{APACrefauthors}%
Ahmadi-Moshkenani, P.%
, Johansen, T\BPBI A.%
\BCBL {}\ \BBA {} Olaru, S.%
\end{APACrefauthors}%
\unskip\
\newblock
\APACrefYearMonthDay{2018}{}{}.
\newblock
{\BBOQ}\APACrefatitle {Combinatorial approach toward multiparametric quadratic
  programming based on characterizing adjacent critical regions} {Combinatorial
  approach toward multiparametric quadratic programming based on characterizing
  adjacent critical regions}.{\BBCQ}
\newblock
\APACjournalVolNumPages{IEEE Transactions on Automatic
  Control}{63}{10}{3221--3231}.
\PrintBackRefs{\CurrentBib}

\bibitem [\protect \citeauthoryear {%
Ahuja%
\ \BBA {} Orlin%
}{%
Ahuja%
\ \BBA {} Orlin%
}{%
{\protect \APACyear {2001}}%
}]{%
Ahuja01}
\APACinsertmetastar {%
Ahuja01}%
\begin{APACrefauthors}%
Ahuja, R\BPBI K.%
\BCBT {}\ \BBA {} Orlin, J\BPBI B.%
\end{APACrefauthors}%
\unskip\
\newblock
\APACrefYearMonthDay{2001}{}{}.
\newblock
{\BBOQ}\APACrefatitle {Inverse optimization} {Inverse optimization}.{\BBCQ}
\newblock
\APACjournalVolNumPages{Operations Research}{49}{5}{771--783}.
\PrintBackRefs{\CurrentBib}

\bibitem [\protect \citeauthoryear {%
Aswani%
, Kaminsky%
, Mintz%
, Flowers%
\BCBL {}\ \BBA {} Fukuoka%
}{%
Aswani%
\ \protect \BOthers {.}}{%
{\protect \APACyear {2019}}%
}]{%
aswani2019behavioral}
\APACinsertmetastar {%
aswani2019behavioral}%
\begin{APACrefauthors}%
Aswani, A.%
, Kaminsky, P.%
, Mintz, Y.%
, Flowers, E.%
\BCBL {}\ \BBA {} Fukuoka, Y.%
\end{APACrefauthors}%
\unskip\
\newblock
\APACrefYearMonthDay{2019}{}{}.
\newblock
{\BBOQ}\APACrefatitle {Behavioral Modeling in Weight Loss Interventions}
  {Behavioral modeling in weight loss interventions}.{\BBCQ}
\newblock
\APACjournalVolNumPages{European Journal of Operational
  Research}{272}{3}{1058--1072}.
\PrintBackRefs{\CurrentBib}

\bibitem [\protect \citeauthoryear {%
Aswani%
, Shen%
\BCBL {}\ \BBA {} Siddiq%
}{%
Aswani%
\ \protect \BOthers {.}}{%
{\protect \APACyear {2018}}%
}]{%
aswani2018inverse}
\APACinsertmetastar {%
aswani2018inverse}%
\begin{APACrefauthors}%
Aswani, A.%
, Shen, Z\BHBI J.%
\BCBL {}\ \BBA {} Siddiq, A.%
\end{APACrefauthors}%
\unskip\
\newblock
\APACrefYearMonthDay{2018}{}{}.
\newblock
{\BBOQ}\APACrefatitle {Inverse optimization with noisy data} {Inverse
  optimization with noisy data}.{\BBCQ}
\newblock
\APACjournalVolNumPages{Operations Research}{66}{3}{870--892}.
\PrintBackRefs{\CurrentBib}

\bibitem [\protect \citeauthoryear {%
Babier%
, Chan%
, Lee%
, Mahmood%
\BCBL {}\ \BBA {} Terekhov%
}{%
Babier%
\ \protect \BOthers {.}}{%
{\protect \APACyear {2021}}%
}]{%
babier2019}
\APACinsertmetastar {%
babier2019}%
\begin{APACrefauthors}%
Babier, A.%
, Chan, T\BPBI C\BPBI Y.%
, Lee, T.%
, Mahmood, R.%
\BCBL {}\ \BBA {} Terekhov, D.%
\end{APACrefauthors}%
\unskip\
\newblock
\APACrefYearMonthDay{2021}{}{}.
\newblock
{\BBOQ}\APACrefatitle {An Ensemble Learning Framework for Model Fitting and
  Evaluation in Inverse Linear Optimization} {An ensemble learning framework
  for model fitting and evaluation in inverse linear optimization}.{\BBCQ}
\newblock
\APACjournalVolNumPages{INFORMS Journal on Optimization}{Articles in
  Advance}{}{1--20}.
\newblock
\begin{APACrefURL} \url{https://doi.org/10.1287/ijoo.2019.0045}
  \end{APACrefURL}
\PrintBackRefs{\CurrentBib}

\bibitem [\protect \citeauthoryear {%
B{\"a}rmann%
, Pokutta%
\BCBL {}\ \BBA {} Schneider%
}{%
B{\"a}rmann%
\ \protect \BOthers {.}}{%
{\protect \APACyear {2017}}%
}]{%
Barmann17}
\APACinsertmetastar {%
Barmann17}%
\begin{APACrefauthors}%
B{\"a}rmann, A.%
, Pokutta, S.%
\BCBL {}\ \BBA {} Schneider, O.%
\end{APACrefauthors}%
\unskip\
\newblock
\APACrefYearMonthDay{2017}{}{}.
\newblock
{\BBOQ}\APACrefatitle {Emulating the Expert: Inverse Optimization through
  Online Learning} {Emulating the expert: Inverse optimization through online
  learning}.{\BBCQ}
\newblock
\BIn{} \APACrefbtitle {International Conference on Machine Learning}
  {International conference on machine learning}\ (\BPGS\ 400--410).
\PrintBackRefs{\CurrentBib}

\bibitem [\protect \citeauthoryear {%
Breiman%
}{%
Breiman%
}{%
{\protect \APACyear {2001}}%
}]{%
breiman2001random}
\APACinsertmetastar {%
breiman2001random}%
\begin{APACrefauthors}%
Breiman, L.%
\end{APACrefauthors}%
\unskip\
\newblock
\APACrefYearMonthDay{2001}{}{}.
\newblock
{\BBOQ}\APACrefatitle {Random forests} {Random forests}.{\BBCQ}
\newblock
\APACjournalVolNumPages{Machine Learning}{45}{1}{5--32}.
\PrintBackRefs{\CurrentBib}

\bibitem [\protect \citeauthoryear {%
Bunduchi%
\ \BBA {} Mandric%
}{%
Bunduchi%
\ \BBA {} Mandric%
}{%
{\protect \APACyear {2011}}%
}]{%
Mathematica}
\APACinsertmetastar {%
Mathematica}%
\begin{APACrefauthors}%
Bunduchi, E.%
\BCBT {}\ \BBA {} Mandric, I.%
\end{APACrefauthors}%
\unskip\
\newblock
\APACrefYearMonthDay{2011}{}{}.
\newblock
\APACrefbtitle {Parametric Linear Programming.} {Parametric linear
  programming.}
\newblock
\begin{APACrefURL}
  \url{http://demonstrations.wolfram.com/ParametricLinearProgramming/}
  \end{APACrefURL}
\PrintBackRefs{\CurrentBib}

\bibitem [\protect \citeauthoryear {%
Burton%
\ \BBA {} Toint%
}{%
Burton%
\ \BBA {} Toint%
}{%
{\protect \APACyear {1992}}%
}]{%
burton1992instance}
\APACinsertmetastar {%
burton1992instance}%
\begin{APACrefauthors}%
Burton, D.%
\BCBT {}\ \BBA {} Toint, P\BPBI L.%
\end{APACrefauthors}%
\unskip\
\newblock
\APACrefYearMonthDay{1992}{}{}.
\newblock
{\BBOQ}\APACrefatitle {On an instance of the inverse shortest paths problem}
  {On an instance of the inverse shortest paths problem}.{\BBCQ}
\newblock
\APACjournalVolNumPages{Mathematical Programming}{53}{1-3}{45--61}.
\PrintBackRefs{\CurrentBib}

\bibitem [\protect \citeauthoryear {%
Burton%
\ \BBA {} Toint%
}{%
Burton%
\ \BBA {} Toint%
}{%
{\protect \APACyear {1994}}%
}]{%
burton_correlated94}
\APACinsertmetastar {%
burton_correlated94}%
\begin{APACrefauthors}%
Burton, D.%
\BCBT {}\ \BBA {} Toint, P\BPBI L.%
\end{APACrefauthors}%
\unskip\
\newblock
\APACrefYearMonthDay{1994}{}{}.
\newblock
{\BBOQ}\APACrefatitle {On the use of an inverse shortest paths algorithm for
  recovering linearly correlated costs} {On the use of an inverse shortest
  paths algorithm for recovering linearly correlated costs}.{\BBCQ}
\newblock
\APACjournalVolNumPages{Mathematical Programming}{63}{1-3}{1--22}.
\PrintBackRefs{\CurrentBib}

\bibitem [\protect \citeauthoryear {%
Chan%
, Craig%
, Lee%
\BCBL {}\ \BBA {} Sharpe%
}{%
Chan%
\ \protect \BOthers {.}}{%
{\protect \APACyear {2014}}%
}]{%
chan2014generalized}
\APACinsertmetastar {%
chan2014generalized}%
\begin{APACrefauthors}%
Chan, T\BPBI C.%
, Craig, T.%
, Lee, T.%
\BCBL {}\ \BBA {} Sharpe, M\BPBI B.%
\end{APACrefauthors}%
\unskip\
\newblock
\APACrefYearMonthDay{2014}{}{}.
\newblock
{\BBOQ}\APACrefatitle {Generalized inverse multiobjective optimization with
  application to cancer therapy} {Generalized inverse multiobjective
  optimization with application to cancer therapy}.{\BBCQ}
\newblock
\APACjournalVolNumPages{Operations Research}{62}{3}{680--695}.
\PrintBackRefs{\CurrentBib}

\bibitem [\protect \citeauthoryear {%
Chan%
, Lee%
\BCBL {}\ \BBA {} Terekhov%
}{%
Chan%
\ \protect \BOthers {.}}{%
{\protect \APACyear {2019}}%
}]{%
GIO}
\APACinsertmetastar {%
GIO}%
\begin{APACrefauthors}%
Chan, T\BPBI C.%
, Lee, T.%
\BCBL {}\ \BBA {} Terekhov, D.%
\end{APACrefauthors}%
\unskip\
\newblock
\APACrefYearMonthDay{2019}{}{}.
\newblock
{\BBOQ}\APACrefatitle {Inverse optimization: Closed-form solutions, geometry,
  and goodness of fit} {Inverse optimization: Closed-form solutions, geometry,
  and goodness of fit}.{\BBCQ}
\newblock
\APACjournalVolNumPages{Management Science}{65}{3}{1115--1135}.
\PrintBackRefs{\CurrentBib}

\bibitem [\protect \citeauthoryear {%
Chow%
\ \BBA {} Recker%
}{%
Chow%
\ \BBA {} Recker%
}{%
{\protect \APACyear {2012}}%
}]{%
chow2012inverse}
\APACinsertmetastar {%
chow2012inverse}%
\begin{APACrefauthors}%
Chow, J.%
\BCBT {}\ \BBA {} Recker, W.%
\end{APACrefauthors}%
\unskip\
\newblock
\APACrefYearMonthDay{2012}{}{}.
\newblock
{\BBOQ}\APACrefatitle {Inverse optimization with endogenous arrival time
  constraints to calibrate the household activity pattern problem} {Inverse
  optimization with endogenous arrival time constraints to calibrate the
  household activity pattern problem}.{\BBCQ}
\newblock
\APACjournalVolNumPages{Transportation Research Part B:
  Methodological}{46}{3}{463--479}.
\PrintBackRefs{\CurrentBib}

\bibitem [\protect \citeauthoryear {%
Drucker%
, Burges%
, Kaufman%
, Smola%
\BCBL {}\ \BBA {} Vapnik%
}{%
Drucker%
\ \protect \BOthers {.}}{%
{\protect \APACyear {1997}}%
}]{%
vapnik95}
\APACinsertmetastar {%
vapnik95}%
\begin{APACrefauthors}%
Drucker, H.%
, Burges, C.%
, Kaufman, L.%
, Smola, A.%
\BCBL {}\ \BBA {} Vapnik, V.%
\end{APACrefauthors}%
\unskip\
\newblock
\APACrefYearMonthDay{1997}{}{}.
\newblock
{\BBOQ}\APACrefatitle {Support vector regression machines} {Support vector
  regression machines}.{\BBCQ}
\newblock
\BIn{} \APACrefbtitle {Advances in Neural Information Processing Systems}
  {Advances in neural information processing systems}\ (\BPGS\ 155--161).
\PrintBackRefs{\CurrentBib}

\bibitem [\protect \citeauthoryear {%
Egri%
, Kis%
, Kov{\'a}cs%
\BCBL {}\ \BBA {} V{\'a}ncza%
}{%
Egri%
\ \protect \BOthers {.}}{%
{\protect \APACyear {2014}}%
}]{%
egri2014inverse}
\APACinsertmetastar {%
egri2014inverse}%
\begin{APACrefauthors}%
Egri, P.%
, Kis, T.%
, Kov{\'a}cs, A.%
\BCBL {}\ \BBA {} V{\'a}ncza, J.%
\end{APACrefauthors}%
\unskip\
\newblock
\APACrefYearMonthDay{2014}{}{}.
\newblock
{\BBOQ}\APACrefatitle {An inverse economic lot-sizing approach to eliciting
  supplier cost parameters} {An inverse economic lot-sizing approach to
  eliciting supplier cost parameters}.{\BBCQ}
\newblock
\APACjournalVolNumPages{International Journal of Production
  Economics}{149}{}{80--88}.
\PrintBackRefs{\CurrentBib}

\bibitem [\protect \citeauthoryear {%
Esfahani%
, Shafieezadeh-Abadeh%
, Hanasusanto%
\BCBL {}\ \BBA {} Kuhn%
}{%
Esfahani%
\ \protect \BOthers {.}}{%
{\protect \APACyear {2018}}%
}]{%
esfahani2018data}
\APACinsertmetastar {%
esfahani2018data}%
\begin{APACrefauthors}%
Esfahani, P.%
, Shafieezadeh-Abadeh, S.%
, Hanasusanto, G.%
\BCBL {}\ \BBA {} Kuhn, D.%
\end{APACrefauthors}%
\unskip\
\newblock
\APACrefYearMonthDay{2018}{}{}.
\newblock
{\BBOQ}\APACrefatitle {Data-driven inverse optimization with imperfect
  information} {Data-driven inverse optimization with imperfect
  information}.{\BBCQ}
\newblock
\APACjournalVolNumPages{Mathematical Programming}{167}{1}{191--234}.
\PrintBackRefs{\CurrentBib}

\bibitem [\protect \citeauthoryear {%
Fern{\'a}ndez-Blanco%
, Morales%
, Pineda%
\BCBL {}\ \BBA {} Porras%
}{%
Fern{\'a}ndez-Blanco%
\ \protect \BOthers {.}}{%
{\protect \APACyear {2019}}%
}]{%
fernandez2019ev}
\APACinsertmetastar {%
fernandez2019ev}%
\begin{APACrefauthors}%
Fern{\'a}ndez-Blanco, R.%
, Morales, J\BPBI M.%
, Pineda, S.%
\BCBL {}\ \BBA {} Porras, {\'A}.%
\end{APACrefauthors}%
\unskip\
\newblock
\APACrefYearMonthDay{2019}{}{}.
\newblock
{\BBOQ}\APACrefatitle {EV-Fleet Power Forecasting via Kernel-Based Inverse
  Optimization} {Ev-fleet power forecasting via kernel-based inverse
  optimization}.{\BBCQ}
\newblock
\APACjournalVolNumPages{arXiv preprint arXiv:1908.00399}{}{}{}.
\PrintBackRefs{\CurrentBib}

\bibitem [\protect \citeauthoryear {%
Goodfellow%
, Bengio%
\BCBL {}\ \BBA {} Courville%
}{%
Goodfellow%
\ \protect \BOthers {.}}{%
{\protect \APACyear {2016}}%
}]{%
Goodfellow-et-al-2016}
\APACinsertmetastar {%
Goodfellow-et-al-2016}%
\begin{APACrefauthors}%
Goodfellow, I.%
, Bengio, Y.%
\BCBL {}\ \BBA {} Courville, A.%
\end{APACrefauthors}%
\unskip\
\newblock
\APACrefYear{2016}.
\newblock
\APACrefbtitle {Deep Learning} {Deep learning}.
\newblock
\APACaddressPublisher{}{MIT Press}.
\newblock
\APACrefnote{\url{http://www.deeplearningbook.org}}
\PrintBackRefs{\CurrentBib}

\bibitem [\protect \citeauthoryear {%
Keshavarz%
, Wang%
\BCBL {}\ \BBA {} Boyd%
}{%
Keshavarz%
\ \protect \BOthers {.}}{%
{\protect \APACyear {2011}}%
}]{%
Keshavarz11}
\APACinsertmetastar {%
Keshavarz11}%
\begin{APACrefauthors}%
Keshavarz, A.%
, Wang, Y.%
\BCBL {}\ \BBA {} Boyd, S.%
\end{APACrefauthors}%
\unskip\
\newblock
\APACrefYearMonthDay{2011}{}{}.
\newblock
{\BBOQ}\APACrefatitle {Imputing a convex objective function} {Imputing a convex
  objective function}.{\BBCQ}
\newblock
\BIn{} \APACrefbtitle {Intelligent Control (ISIC), 2011 IEEE International
  Symposium on} {Intelligent control (isic), 2011 ieee international symposium
  on}\ (\BPGS\ 613--619).
\PrintBackRefs{\CurrentBib}

\bibitem [\protect \citeauthoryear {%
Kov{\'a}cs%
}{%
Kov{\'a}cs%
}{%
{\protect \APACyear {2019}}%
}]{%
kovacs2019parameter}
\APACinsertmetastar {%
kovacs2019parameter}%
\begin{APACrefauthors}%
Kov{\'a}cs, A.%
\end{APACrefauthors}%
\unskip\
\newblock
\APACrefYearMonthDay{2019}{}{}.
\newblock
{\BBOQ}\APACrefatitle {Parameter Elicitation for Consumer Models in Demand
  Response Management} {Parameter elicitation for consumer models in demand
  response management}.{\BBCQ}
\newblock
\BIn{} \APACrefbtitle {2019 1st Global Power, Energy and Communication
  Conference (GPECOM)} {2019 1st global power, energy and communication
  conference (gpecom)}\ (\BPGS\ 456--460).
\PrintBackRefs{\CurrentBib}

\bibitem [\protect \citeauthoryear {%
Oberdieck%
, Diangelakis%
, Papathanasiou%
, Nascu%
\BCBL {}\ \BBA {} Pistikopoulos%
}{%
Oberdieck%
\ \protect \BOthers {.}}{%
{\protect \APACyear {2016}}%
}]{%
POP}
\APACinsertmetastar {%
POP}%
\begin{APACrefauthors}%
Oberdieck, R.%
, Diangelakis, N.%
, Papathanasiou, M.%
, Nascu, I.%
\BCBL {}\ \BBA {} Pistikopoulos, E.%
\end{APACrefauthors}%
\unskip\
\newblock
\APACrefYearMonthDay{2016}{}{}.
\newblock
{\BBOQ}\APACrefatitle {{POP --} {P}arametric {O}ptimization {T}oolbox} {{POP
  --} {P}arametric {O}ptimization {T}oolbox}.{\BBCQ}
\newblock
\APACjournalVolNumPages{Industrial \& Engineering Chemistry
  Research}{55}{33}{8979-8991}.
\newblock
\begin{APACrefURL} \url{https://doi.org/10.1021/acs.iecr.6b01913}
  \end{APACrefURL}
\newblock
\begin{APACrefDOI} \doi{10.1021/acs.iecr.6b01913} \end{APACrefDOI}
\PrintBackRefs{\CurrentBib}

\bibitem [\protect \citeauthoryear {%
Pedregosa%
\ \protect \BOthers {.}}{%
Pedregosa%
\ \protect \BOthers {.}}{%
{\protect \APACyear {2011}}%
}]{%
scikit}
\APACinsertmetastar {%
scikit}%
\begin{APACrefauthors}%
Pedregosa, F.%
, Varoquaux, G.%
, Gramfort, A.%
, Michel, V.%
, Thirion, B.%
, Grisel, O.%
\BDBL {}Duchesnay, E.%
\end{APACrefauthors}%
\unskip\
\newblock
\APACrefYearMonthDay{2011}{{\APACmonth{11}}}{}.
\newblock
{\BBOQ}\APACrefatitle {Scikit-Learn: Machine Learning in Python} {Scikit-learn:
  Machine learning in python}.{\BBCQ}
\newblock
\APACjournalVolNumPages{J. Mach. Learn. Res.}{12}{null}{2825–2830}.
\PrintBackRefs{\CurrentBib}

\bibitem [\protect \citeauthoryear {%
Pistikopoulos%
, Georgiadis%
\BCBL {}\ \BBA {} Dua%
}{%
Pistikopoulos%
\ \protect \BOthers {.}}{%
{\protect \APACyear {2011}}%
}]{%
pistikopoulos2011multi}
\APACinsertmetastar {%
pistikopoulos2011multi}%
\begin{APACrefauthors}%
Pistikopoulos, E\BPBI N.%
, Georgiadis, M\BPBI C.%
\BCBL {}\ \BBA {} Dua, V.%
\end{APACrefauthors}%
\unskip\
\newblock
\APACrefYear{2011}.
\newblock
\APACrefbtitle {Multi-Parametric Programming} {Multi-parametric programming}\
  (\BVOL~1).
\PrintBackRefs{\CurrentBib}

\bibitem [\protect \citeauthoryear {%
Rasmussen%
}{%
Rasmussen%
}{%
{\protect \APACyear {2004}}%
}]{%
rasmussen2004gaussian}
\APACinsertmetastar {%
rasmussen2004gaussian}%
\begin{APACrefauthors}%
Rasmussen, C.%
\end{APACrefauthors}%
\unskip\
\newblock
\APACrefYearMonthDay{2004}{}{}.
\newblock
{\BBOQ}\APACrefatitle {Gaussian processes in machine learning} {Gaussian
  processes in machine learning}.{\BBCQ}
\newblock
\BIn{} \APACrefbtitle {Advanced Lectures on Machine Learning} {Advanced
  lectures on machine learning}\ (\BPGS\ 63--71).
\newblock
\APACaddressPublisher{}{Springer}.
\PrintBackRefs{\CurrentBib}

\bibitem [\protect \citeauthoryear {%
Russell%
\ \BBA {} Norvig%
}{%
Russell%
\ \BBA {} Norvig%
}{%
{\protect \APACyear {2016}}%
}]{%
russell2016artificial}
\APACinsertmetastar {%
russell2016artificial}%
\begin{APACrefauthors}%
Russell, S\BPBI J.%
\BCBT {}\ \BBA {} Norvig, P.%
\end{APACrefauthors}%
\unskip\
\newblock
\APACrefYear{2016}.
\newblock
\APACrefbtitle {Artificial intelligence: a modern approach} {Artificial
  intelligence: a modern approach}.
\newblock
\APACaddressPublisher{}{Malaysia; Pearson Education Limited,}.
\PrintBackRefs{\CurrentBib}

\bibitem [\protect \citeauthoryear {%
Saez-Gallego%
\ \BBA {} Morales%
}{%
Saez-Gallego%
\ \BBA {} Morales%
}{%
{\protect \APACyear {2017}}%
}]{%
Gallego}
\APACinsertmetastar {%
Gallego}%
\begin{APACrefauthors}%
Saez-Gallego, J.%
\BCBT {}\ \BBA {} Morales, J\BPBI M.%
\end{APACrefauthors}%
\unskip\
\newblock
\APACrefYearMonthDay{2017}{}{}.
\newblock
{\BBOQ}\APACrefatitle {Short-term forecasting of price-responsive loads using
  inverse optimization} {Short-term forecasting of price-responsive loads using
  inverse optimization}.{\BBCQ}
\newblock
\APACjournalVolNumPages{{IEEE} Transactions on Smart Grid}{9}{5}{4805--4814}.
\PrintBackRefs{\CurrentBib}

\bibitem [\protect \citeauthoryear {%
Shahmoradi%
\ \BBA {} Lee%
}{%
Shahmoradi%
\ \BBA {} Lee%
}{%
{\protect \APACyear {2020}}%
}]{%
shahmoradi2020}
\APACinsertmetastar {%
shahmoradi2020}%
\begin{APACrefauthors}%
Shahmoradi, Z.%
\BCBT {}\ \BBA {} Lee, T.%
\end{APACrefauthors}%
\unskip\
\newblock
\APACrefYearMonthDay{2020}{}{}.
\newblock
\APACrefbtitle {Quantile Inverse Optimization: Improving Stability in Inverse
  Linear Programming.} {Quantile inverse optimization: Improving stability in
  inverse linear programming.}
\PrintBackRefs{\CurrentBib}

\bibitem [\protect \citeauthoryear {%
Tan%
, Delong%
\BCBL {}\ \BBA {} Terekhov%
}{%
Tan%
\ \protect \BOthers {.}}{%
{\protect \APACyear {2019}}%
}]{%
deepIO}
\APACinsertmetastar {%
deepIO}%
\begin{APACrefauthors}%
Tan, Y.%
, Delong, A.%
\BCBL {}\ \BBA {} Terekhov, D.%
\end{APACrefauthors}%
\unskip\
\newblock
\APACrefYearMonthDay{2019}{}{}.
\newblock
{\BBOQ}\APACrefatitle {Deep Inverse Optimization} {Deep inverse
  optimization}.{\BBCQ}
\newblock
\BIn{} \APACrefbtitle {International Conference on Integration of Constraint
  Programming, Artificial Intelligence, and Operations Research} {International
  conference on integration of constraint programming, artificial intelligence,
  and operations research}\ (\BPGS\ 540--556).
\PrintBackRefs{\CurrentBib}

\bibitem [\protect \citeauthoryear {%
Tan%
, Terekhov%
\BCBL {}\ \BBA {} Delong%
}{%
Tan%
\ \protect \BOthers {.}}{%
{\protect \APACyear {2020}}%
}]{%
tan2020}
\APACinsertmetastar {%
tan2020}%
\begin{APACrefauthors}%
Tan, Y.%
, Terekhov, D.%
\BCBL {}\ \BBA {} Delong, A.%
\end{APACrefauthors}%
\unskip\
\newblock
\APACrefYearMonthDay{2020}{}{}.
\newblock
{\BBOQ}\APACrefatitle {Learning linear programs from optimal decisions}
  {Learning linear programs from optimal decisions}.{\BBCQ}
\newblock
\BIn{} \APACrefbtitle {Advances in Neural Information Processing Systems}
  {Advances in neural information processing systems}\ (\BVOL~33).
\PrintBackRefs{\CurrentBib}

\bibitem [\protect \citeauthoryear {%
Taskar%
, Chatalbashev%
, Koller%
\BCBL {}\ \BBA {} Guestrin%
}{%
Taskar%
\ \protect \BOthers {.}}{%
{\protect \APACyear {2005}}%
}]{%
taskar2005learning}
\APACinsertmetastar {%
taskar2005learning}%
\begin{APACrefauthors}%
Taskar, B.%
, Chatalbashev, V.%
, Koller, D.%
\BCBL {}\ \BBA {} Guestrin, C.%
\end{APACrefauthors}%
\unskip\
\newblock
\APACrefYearMonthDay{2005}{}{}.
\newblock
{\BBOQ}\APACrefatitle {Learning Structured Prediction Models: A Large Margin
  Approach} {Learning structured prediction models: A large margin
  approach}.{\BBCQ}
\newblock
\BIn{} \APACrefbtitle {Proceedings of the 22nd International Conference on
  Machine Learning} {Proceedings of the 22nd international conference on
  machine learning}\ (\BPGS\ 896--903).
\PrintBackRefs{\CurrentBib}

\bibitem [\protect \citeauthoryear {%
Tavasl{\i}o{\u{g}}lu%
, Lee%
, Valeva%
\BCBL {}\ \BBA {} Schaefer%
}{%
Tavasl{\i}o{\u{g}}lu%
\ \protect \BOthers {.}}{%
{\protect \APACyear {2018}}%
}]{%
tavasliouglu2018}
\APACinsertmetastar {%
tavasliouglu2018}%
\begin{APACrefauthors}%
Tavasl{\i}o{\u{g}}lu, O.%
, Lee, T.%
, Valeva, S.%
\BCBL {}\ \BBA {} Schaefer, A\BPBI J.%
\end{APACrefauthors}%
\unskip\
\newblock
\APACrefYearMonthDay{2018}{}{}.
\newblock
{\BBOQ}\APACrefatitle {On the structure of the inverse-feasible region of a
  linear program} {On the structure of the inverse-feasible region of a linear
  program}.{\BBCQ}
\newblock
\APACjournalVolNumPages{Operations Research Letters}{46}{1}{147--152}.
\PrintBackRefs{\CurrentBib}

\end{thebibliography}

\end{document}